# A NEW METHOD OF NORMAL APPROXIMATION[1]

### By Sourav Chatterjee

#### *University of California, Berkeley*

We introduce a new version of Stein's method that reduces a large class of normal approximation problems to variance bounding exercises, thus making a connection between central limit theorems and concentration of measure. Unlike Skorokhod embeddings, the object whose variance must be bounded has an explicit formula that makes it possible to carry out the program more easily. As an application, we derive a general CLT for functions that are obtained as combinations of many local contributions, where the definition of "local" itself depends on the data. Several examples are given, including the solution to a nearest-neighbor CLT problem posed by P. Bickel.

**1. Introduction.** Central limit theorems for general nonadditive functions of independent random variables have been studied by various authors using a variety of techniques. Some examples are: (i) the method of Hajék projections and some sophisticated extensions (e.g., [20, 38, 43]); (ii) Stein's method of normal approximation (references in Section 2); (iii) the big-blocks–small-blocks technique and its modern multidimensional versions (e.g. [2, 6]); (iv) the martingale approach and Skorokhod embeddings; (v) the method of moments. In this paper, we present a new approach that may go beyond the limitations of these existing techniques. The power of the method is demonstrated through several applications, mainly geometrical in nature, that are otherwise difficult. In the related article [10], we provide some applications to random matrices.

The paper is organized as follows. Section 2 contains a brief discussion of Stein's method and our main results (Theorems 2.2 and 2.5). Examples are worked out in Section 3. Proofs of the main theorems are in Section 4.

Received November 2006; revised July 2007.

[1]Supported in part by NSF Grant DMS-07-07054 and a Sloan Research Fellowship in Mathematics.

**2. Results.** Recall that the Kantorovich–Wasserstein distance between two probability measures $\mu$ and $\nu$ on the real line is defined as

$$\mathcal{W}(\mu, \nu) := \sup \left\{ \left| \int h \, d\mu - \int h \, d\nu \right| : h \text{ Lipschitz, with } \|h\|_{\text{Lip}} \leq 1 \right\}.$$

Convergence of measures in the Kantorovich–Wasserstein metric is stronger than weak convergence. Based on the Kantorovich–Wasserstein distance, we introduce the following measure of "distance to Gaussianity."

DEFINITION 2.1. Let $W$ be a real-valued random variable with finite second moment. Let $\mu$ be the law of $(W - \mathbb{E}(W))/\sqrt{\text{Var}(W)}$ and let $\nu$ be the standard Gaussian law. We define

$$\delta_W := \mathcal{W}(\mu, \nu).$$

This is our preferred metric of normal approximation in this paper. Using the bounds on the Wasserstein distance, analogous results can be obtained for the Kolmogorov distance via smoothing, but the rates will be suboptimal. This problem is very common in Stein's method; obtaining optimal rates for the Kolmogorov metric requires extra work and new ideas (see, e.g., [13]). Since our main focus is on convergence to normality and not so much on error bounds, we will not worry about this issue here.

2.1. *Stein's method.* A well-known computation via integration by parts shows that if $Z \sim N(0, 1)$, then $\mathbb{E}(\varphi(Z)Z) = \mathbb{E}(\varphi'(Z))$ for all absolutely continuous $\varphi$ with $\mathbb{E}|\varphi'(Z)| < \infty$. Conversely, if $W$ is a random variable satisfying $\mathbb{E}(\varphi(W)W) = \mathbb{E}(\varphi'(W))$ for all Lipschitz $\varphi$, then $W \sim N(0, 1)$. Consequently, if $W$ is such that $\mathbb{E}(\varphi(W)W) \approx \mathbb{E}(\varphi'(W))$ for all $\varphi$ belonging to a large class of functions, then one can expect the distribution of $W$ to be close to the $N(0, 1)$ distribution.

This is the key idea behind Stein's method of normal approximation, introduced by Charles Stein in the seminal paper [40] and later developed in his book [41]. Precise error bounds can be obtained in various ways. We will reproduce one of Stein's results (Lemma 4.2) that gives a bound on the Kantorovich–Wasserstein distance to normality.

However, the problem begins at this point. Given a random variable $W$ that may be a complicated function of many other variables, there is no general method for showing that $\mathbb{E}(\varphi(W)W) \approx \mathbb{E}(\varphi'(W))$. Several powerful techniques for carrying out this step under special conditions are available in the literature on Stein's method (e.g., exchangeable pairs [41], diffusion generators [5], dependency graphs [4, 13, 36], size bias transforms [23], zero bias transforms [22], specialized procedures like [21, 34, 35] and recent developments [11, 12], to cite a few), but, somehow, they all require something

"nice" to happen. It is rarely the case that something arbitrary (like the Levina–Bickel statistic [30], to be discussed later) becomes amenable to any of the existing versions of Stein's method.

This is the ground that we attempt to break in this paper. Given a random variable $W$ that is an explicit but arbitrary function of a collection of independent random variables, satisfying $\mathbb{E}(W) = 0$ and $\mathbb{E}(W^2) = 1$, we prescribe a method of constructing another random variable $T$ so that for all smooth $\varphi$, we have

$$\mathbb{E}(\varphi(W)W) \approx \mathbb{E}(\varphi'(W)T).$$

In particular, taking $\varphi$ to be the identity function, we get $\mathbb{E}(T) \approx 1$. If, now, $\mathrm{Var}(T)$ is small enough, then we can make the easy but crucial deduction that $T$ "can be substituted by the constant $\mathbb{E}(T)$" to get $\mathbb{E}(\varphi(W)W) \approx \mathbb{E}(\varphi'(W))$, which shows that the distribution of $W$ is approximately standard Gaussian. Thus, the normal approximation problem is reduced to the problem of bounding the variance of $T$. Of course, the crux of the matter lies in the construction of $T$, which we undertake below.

2.2. *An abstract result.* Let $\mathcal{X}$ be a measure space and suppose $X = (X_1, \ldots, X_n)$ is a vector of independent $\mathcal{X}$-valued random variables. Let $X' = (X_1', \ldots, X_n')$ be an independent copy of $X$. Let $[n] = \{1, \ldots, n\}$, and for each $A \subseteq [n]$, define the random vector $X^A$ as

$$X_i^A = \begin{cases} X_i', & \text{if } i \in A, \\ X_i, & \text{if } i \notin A. \end{cases}$$

When $A$ is a singleton set like $\{j\}$, we will simply write $X^j$ instead of $X^{\{j\}}$. Let $f \colon \mathcal{X} \to \mathbb{R}$ be a measurable function. We define a randomized derivative of $f(X)$ along the $j$th coordinate as

$$\Delta_j f(X) := f(X) - f(X^j).$$

Note that $\Delta_j f(X)$ depends not only on the vector $X$, but also on $X_j'$. Next, for each $A \subseteq [n]$, let

$$T_A := \sum_{j \notin A} \Delta_j f(X) \Delta_j f(X^A)$$

and let

(1) $$T = \frac{1}{2} \sum_{A \subsetneq [n]} \frac{T_A}{\binom{n}{|A|}(n - |A|)}.$$

Putting $W = f(X)$ and assuming that $\mathbb{E}(W) = 0$, we show (Lemma 2.4) that whenever $\sum_j \mathbb{E}|\Delta_j f(X)|^3$ is small, we have $\mathbb{E}(\varphi(W)W) \approx \mathbb{E}(\varphi'(W)T)$ for all $\varphi$ belonging to a large class of functions. The main consequence is the following normal approximation theorem, which is our main abstract result.

THEOREM 2.2. *Let all terms be defined as above and let $W = f(X)$. Suppose that $\mathbb{E}(W) = 0$ and $\sigma^2 := \mathbb{E}(W^2) < \infty$. Then, $\mathbb{E}(T) = \sigma^2$ and*

$$\delta_W \leq \frac{[\mathrm{Var}(\mathbb{E}(T|W))]^{1/2}}{\sigma^2} + \frac{1}{2\sigma^3} \sum_{j=1}^{n} \mathbb{E}|\Delta_j f(X)|^3,$$

*where $\delta_W$ is the distance to normality defined in Definition 2.1.*

For the simplest possible application of Theorem 2.2, let $X_1, \ldots, X_n$ be i.i.d. real-valued random variables with $\mathbb{E}(X_i) = 0$ and $\mathbb{E}(X_i^2) = 1$, and let $W = f(X) := n^{-1/2} \sum_{i=1}^{n} X_i$. Then, for any $A \subseteq [n]$ and $j \notin A$,

$$\Delta_j f(X^A) = n^{-1/2}(X_j - X_j').$$

Thus, $T_A = n^{-1} \sum_{j \notin A} (X_j - X_j')^2$. A simple verification now shows that

$$T = \frac{1}{2n} \sum_{j=1}^{n} (X_j - X_j')^2.$$

Now, assuming that $\mathbb{E}(X_i^4) < \infty$ and using the inequality $\mathrm{Var}(\mathbb{E}(T|W)) \leq \mathrm{Var}(T)$, we see that $\delta_W \leq Cn^{-1/2}$ for some constant $C$ that depends only on the distribution of the $X_i$'s.

A shortcoming of Theorem 2.2 is that it does not say anything about the variance $\sigma^2$. Somewhat mysteriously, we get a normal approximation result without having to evaluate the variance of our statistic. Of course, the error bound depends on $\sigma^2$ and to show that the bound is useful, we require a lower bound on $\sigma^2$. We prefer to think of that as a separate problem.

The proof of Theorem 2.2, and indeed our whole technique, rests on the following "local-to-global" lemma that deserves to be mentioned its own right. It is closely connected to certain techniques introduced in the author's previous works [8, 9].

LEMMA 2.3. *For any $g$, $f: \mathcal{X}^n \to \mathbb{R}$ such that $\mathbb{E}g(X)^2$ and $\mathbb{E}f(X)^2$ are both finite, we have*

$$\mathrm{Cov}(g(X), f(X)) = \frac{1}{2} \sum_{A \subsetneq [n]} \frac{1}{\binom{n}{|A|}(n - |A|)} \sum_{j \notin A} \mathbb{E}[\Delta_j g(X) \Delta_j f(X^A)].$$

A consequence of the above lemma, to be proven in Section 4, is the following result, which shows that $\mathbb{E}(W\varphi(W)) \approx \mathbb{E}(\varphi'(W)T)$ for all "nice" $\varphi$, whenever $\sum_{i=1}^{n} \mathbb{E}|\Delta_i f(X)|^3$ is small.

LEMMA 2.4. *Let $W = f(X)$ and suppose that $\mathbb{E}(W) = 0$ and $\mathbb{E}(W^2) = 1$. Then, for any $\varphi \in C^2(\mathbb{R})$ with bounded second derivative, we have*

$$|\mathbb{E}(\varphi(W)W) - \mathbb{E}(\varphi'(W)T)| \leq \frac{\|\varphi''\|_\infty}{4} \sum_{j=1}^{n} \mathbb{E}|\Delta_j f(X)|^3.$$

This lemma has a connection with the Goldstein–Reinert zero-bias transform version of Stein's method [22], which we now explain. Given a random variable $W$ with mean zero and unit variance, a random variable $W^*$ is said to be a zero-bias transform of $W$ if for all absolutely continuous $\varphi$,

$$\mathbb{E}(W\varphi(W)) = \mathbb{E}(\varphi'(W^*))$$

whenever both sides are well defined. It is shown in [22] that a zero-bias transform always exists and the closeness to normality for $W$ can be measured by the closeness of the distributions of $W$ and $W^*$ (which is usually done by constructing $W^*$ such that $W \approx W^*$). The problem with this approach, again, is that zero-bias transforms are hard to construct in general. Lemma 2.4 tells us that whenever $\sum_{i=1}^n \mathbb{E}|\Delta_i f(X)|^3$ is small, we have $\mathbb{E}(W\varphi(W)) \approx \mathbb{E}(\varphi'(W)\mathbb{E}(T|W))$. This means, roughly, that $\mathbb{E}(T|W)$ is approximately the Radon–Nikodym density of the law of $W^*$ with respect to the law of $W$, although such a density may not actually exist. Incidentally, such densities have been studied before, for example, in [7].

### 2.3. A general CLT for structures with local dependence.

Numerous central limit theorems in probability theory have been conjectured or proven by following the intuition that a CLT for a sum of dependent summands should hold if "the dependencies are local in nature." Some notable examples are the classical big-blocks–small-blocks technique for analyzing $m$-dependent sequences, its multidimensional generalizations (e.g., [2, 6]), and the dependency graph method of [3]. Here, we provide a new method that is seemingly more powerful than the existing techniques (our applications provide some evidence for this claim) and also gives explicit error bounds. The method is derived as a nontrivial corollary of Theorem 2.2.

Let $\mathcal{X}$ be a measure space and suppose $f : \mathcal{X}^n \to \mathbb{R}$ is a measurable map, where $n \geq 1$ is a fixed positive integer. Suppose $G$ is a map which associates to every $x \in \mathcal{X}^n$ an undirected graph $G(x)$ on $[n] := \{1, \ldots, n\}$. Such a map will be called a *graphical rule* on $\mathcal{X}^n$. We will say that a graphical rule $G$ is *symmetric* if for any permutation $\pi$ of $[n]$ and any $(x_1, \ldots, x_n) \in \mathcal{X}^n$, the set of edges in $G(x_{\pi(1)}, \ldots, x_{\pi(n)})$ is exactly

$$\{\{\pi(i), \pi(j)\} : \{i, j\} \in G(x_1, \ldots, x_n)\}.$$

Now, fix $m > n$. We say that a vector $x \in \mathcal{X}^n$ is *embedded* in another vector $y \in \mathcal{X}^m$ if there exist distinct $i_1, \ldots, i_n \in [m]$ with $x_k = y_{i_k}$ for $1 \leq k \leq n$. A graphical rule $G'$ on $\mathcal{X}^m$ will be called an *extension* of $G$ if for any $x \in \mathcal{X}^n$ embedded in $y \in \mathcal{X}^m$, the graph $G(x)$ on $[n]$ is the naturally induced subgraph of the graph $G'(y)$ on $[m]$.

Now, take any $x, x' \in \mathcal{X}^n$. For each $i \in [n]$, let $x^i$ be the vector obtained by replacing $x_i$ with $x_i'$ in the vector $x$. For any two distinct elements $i$ and $j$ of $[n]$, let $x^{ij}$ be the vector obtained by replacing $x_i$ with $x_i'$ and $x_j$ with $x_j'$. We say that the coordinates $i$ and $j$ are *noninteracting* under the triple $(f, x, x')$ if

$$f(x) - f(x^j) = f(x^i) - f(x^{ij}).$$

Note that the definition is symmetric in $i$ and $j$. This is just a discrete analog of the condition

$$\frac{\partial^2 f}{\partial x_i \, \partial x_j}(x) = 0,$$

which clarifies why it is reasonable to define interaction between coordinates in this manner.

We will say that a graphical rule $G$ is an *interaction rule* for a function $f$ if for any choice of $x, x'$ and $i, j$, the event that $\{i, j\}$ is *not* an edge in the graphs $G(x)$, $G(x^i)$, $G(x^j)$ and $G(x^{ij})$ implies that $i$ and $j$ are noninteracting vertices under $(f, x, x')$. Again, in a continuous setup, we would simply declare that $G(x)$ is the graph that puts an edge between $i$ and $j$ if and only if

$$\frac{\partial^2 f}{\partial x_i \, \partial x_j}(x) \neq 0.$$

Clearly, this is a naturally acceptable definition of an interaction rule (or interaction graph) for $f$. Since we do not want to confine ourselves to the continuous case, the definitions become a bit more complex.

THEOREM 2.5. *Let* $f : \mathcal{X}^n \to \mathbb{R}$ *be a measurable map that admits a symmetric interaction rule* $G$. *Let* $X_1, X_2, \ldots$ *be a sequence of i.i.d.* $\mathcal{X}$-*valued random variables and let* $X = (X_1, \ldots, X_n)$. *Let* $W = f(X)$ *and* $\sigma^2 = \mathrm{Var}(W)$. *Let* $X' = (X_1', \ldots, X_n')$ *be an independent copy of* $X$. *For each* $j$, *define*

$$\Delta_j f(X) = W - f(X_1, \ldots, X_{j-1}, X_j', X_{j+1}, \ldots, X_n)$$

*and let* $M = \max_j |\Delta_j f(X)|$. *Let* $G'$ *be an arbitrary symmetric extension of* $G$ *on* $\mathcal{X}^{n+4}$ *and put*

$$\delta := 1 + \text{degree of the vertex } 1 \text{ in } G'(X_1, \ldots, X_{n+4}).$$

*We then have*

$$\delta_W \leq \frac{Cn^{1/2}}{\sigma^2} \mathbb{E}(M^8)^{1/4} \mathbb{E}(\delta^4)^{1/4} + \frac{1}{2\sigma^3} \sum_{j=1}^n \mathbb{E}|\Delta_j f(X)|^3,$$

*where* $\delta_W$ *is the distance to normality defined in Definition 2.1 and* $C$ *is a universal constant.*

**3. Examples.** This section is devoted to working out applications of Theorems 2.2 and 2.5. Some of these are new results, while others are simpler proofs of existing results. In general, we do not investigate whether our convergence rates are optimal, but in examples where the answers are known, our rates match the existing ones. References to the relevant literature are given in the appropriate places.

3.1. *Quadratic forms.* Suppose $X_1, \ldots, X_n$ are i.i.d. real-valued random variables with zero mean, unit variance and finite fourth moment. Let $\mathbf{A} = (a_{ij})_{1 \leq i, j \leq n}$ be a real symmetric matrix. We consider the following question: under what conditions on the matrix $\mathbf{A}$ can we say that the quadratic form $W = \sum_{i \leq j} a_{ij} X_i X_j$ is approximately Gaussian?

The answer to this question is not very simple; for instance, the usual methods for U-statistics do not work for this problem. The best known condition in the literature (see, e.g., Rotar [37], Hall [27], de Jong [15]) says that asymptotic normality holds if we have a sequence of symmetric matrices $\mathbf{A}_n$ satisfying

$$(2) \qquad \lim_{n \to \infty} \sigma_n^{-4} \operatorname{Tr}(\mathbf{A}_n^4) = 0 \quad \text{and} \quad \lim_{n \to \infty} \sigma_n^{-2} \max_i \sum_{j=1}^n a_{n,ij}^2 = 0,$$

where $\sigma_n^2 = \frac{1}{2} \operatorname{Tr}(\mathbf{A}_n^2) = \operatorname{Var}(W_n)$. The first condition may seem strange, but it is actually equivalent to

$$\lim_{n \to \infty} \frac{\mathbb{E}(W_n - \mathbb{E}(W_n))^4}{(\operatorname{Var}(W_n))^2} = 3,$$

which is a necessary condition for convergence to normality if the sequence $\{W_n^4\}_{n \geq 1}$ is uniformly integrable. The best error bounds were obtained by Götze and Tikhomirov [24, 25].

It is possible to deal with this problem quite easily using our method. Since this is meant to be only an illustration, we keep the expressions as simple as possible by letting the $X_i$'s be $\pm 1$ Rademacher random variables.

PROPOSITION 3.1. *Let $X = (X_1, \ldots, X_n)$ be a vector of i.i.d. random variables with $\mathbb{P}(X_i = 1) = \mathbb{P}(X_i = -1) = 1/2$. Let $\mathbf{A} = (a_{ij})_{1 \leq i, j \leq n}$ be a real symmetric matrix. Let $W = \sum_{i \leq j} a_{ij} X_i X_j$ and $\sigma^2 = \operatorname{Var}(W) = \frac{1}{2} \operatorname{Tr}(\mathbf{A}^2)$. Then,*

$$\delta_W \leq \left( \frac{\operatorname{Tr}(\mathbf{A}^4)}{2\sigma^4} \right)^{1/2} + \frac{5}{2\sigma^3} \sum_{i=1}^n \left( \sum_{j=1}^n a_{ij}^2 \right)^{3/2},$$

*where $\delta_W$ is the distance to normality defined in Definition 2.1.*

Note that the classical condition (2) is implied by the above result, because

$$\sum_{i=1}^n \left(\sum_{j=1}^n a_{ij}^2\right)^{3/2} \le 2\sigma^2 \max_i \left(\sum_{j=1}^n a_{ij}^2\right)^{1/2}.$$

PROOF OF PROPOSITION 3.1.   We will freely use the notation from Theorem 2.2 in this proof. Without loss of generality, we can replace $a_{ij}$ by $a_{ij}/\sigma$ and assume that $\mathrm{Tr}(\mathbf{A}^2) = \sum_{i,j} a_{ij}^2 = 2$. Again, since $\mathbb{E}(W) = \sum_{i=1}^n a_{ii}$, we can assume that $a_{ii} = 0$ for all $i$ after subtracting the mean. Then, note that for any $A \subseteq [n]$ and $i \notin A$,

$$\Delta_i f(X^A) = (X_i - X_i')\left(\sum_{j \notin A} a_{ij}X_j + \sum_{j \in A} a_{ij}X_j'\right).$$

Thus, we have

$$\mathbb{E}(\Delta_i f(X)\Delta_i f(X^A)|X)$$
$$= \mathbb{E}\left((X_i - X_i')^2\left(\sum_{j=1}^n a_{ij}X_j\right)\left(\sum_{j \notin A} a_{ij}X_j + \sum_{j \in A} a_{ij}X_j'\right)\Big|X\right)$$
$$= 2\left(\sum_{j=1}^n a_{ij}X_j\right)\left(\sum_{j \notin A} a_{ij}X_j\right) = 2\sum_{j \in [n]\setminus A, k \in [n]} a_{ij}a_{ik}X_jX_k.$$

A simple verification now shows that

$$\mathbb{E}(T|X) = \sum_{1 \le i,j,k \le n} a_{ij}a_{ik}X_jX_k\left(\sum_{A \subseteq [n]\setminus\{i,j\}} \frac{1}{\binom{n}{|A|}(n - |A|)}\right)$$
$$= \frac{1}{2}X^t\mathbf{A}^2 X,$$

where $X^t$ stands for the transpose of the column vector $X$. Let $b_{ij}$ denote the $(i,j)$th element of $\mathbf{A}^2$. Since $X_i^2 \equiv 1$, the above identity shows that

$$\mathrm{Var}(\mathbb{E}(T|X)) = \mathrm{Var}\left(\sum_{i<j} b_{ij}X_iX_j\right) = \sum_{i<j} b_{ij}^2 \le \tfrac{1}{2}\mathrm{Tr}(\mathbf{A}^4).$$

Finally, by Khintchine's inequality [26], we get

$$\mathbb{E}|\Delta_i f(X)|^3 = 4\mathbb{E}\left|\sum_{j=1}^n a_{ij}X_j\right|^3 \le 5\left(\sum_{j=1}^n a_{ij}^2\right)^{3/2}.$$

The proof is now completed by using the above bounds in Theorem 2.2.   □

3.2. *An occupancy problem.* Suppose $n$ balls are dropped into $\alpha n$ boxes such that all $(\alpha n)^n$ possibilities are equally likely. Let $W$ be the number of empty boxes. The distribution of $W$ is completely known from elementary probability (see, e.g., Feller [19], Section IV.2; for extensive references, see [16]). Very general error bounds for the normal approximation of random variables like $W$ are also known [18]. For illustrative purposes, we now apply Theorem 2.5 to prove a CLT for $W$ when $\alpha$ remains fixed and $n$ tends to infinity.

PROPOSITION 3.2. *Let $W$ be the number of empty boxes as above. Then,*

$$\delta_W \le \frac{C f(\alpha)}{\sqrt{n}},$$

*where $\delta_W$ is the distance to normality defined in Definition 2.1, $f(\alpha) = (\alpha e^{-1/\alpha} - (1 + \alpha)e^{-2/\alpha})^{-3/2}$ and $C$ is a universal constant.*

REMARK. This matches the sharp convergence rate obtained in [18], although that result is for the Kolmogorov distance.

PROOF OF PROPOSITION 3.2. In the following discussion, we are going to freely use the terms defined in the statement of Theorem 2.5 without explicit mention. Let $\mathcal{X}$ be the set of labels of the $\alpha n$ boxes and let $X_i$ denote the label of the box into which ball $i$ is dropped. Let $X = (X_1, \ldots, X_n)$ and let $W = f(X)$ denote the number of empty boxes in the configuration $X$. Then, the transformation $X \to X^j$ denotes the action of moving the ball $j$ from its current box to a box chosen uniformly at random. Clearly, $|\Delta_j f(X)| \le 1$ always and therefore $M \le 1$, where $M = \max_j |\Delta_j f(X)|$ as defined in Theorem 2.5.

Let us now define an interaction graph for this problem. Given a configuration $x$, let $G(x)$ be the graph on $[n]$ that puts an edge between $i$ and $j$ if and only if $x_i = x_j$, that is, the balls $i$ and $j$ land in the same box in the configuration $x$. It is easy to see that $G$ is symmetric. Let us show that $G$ is indeed an interaction graph for $f$ according to our definition.

Let $x'$ be another configuration and let $x^i$, $x^j$ and $x^{ij}$ be defined as usual. Suppose $\{i, j\}$ is not an edge in $G(x)$, $G(x^i)$, $G(x^j)$ and $G(x^{ij})$. This means that the balls $i$ and $j$ are in different boxes in all four configurations. Now, $f(x) - f(x^j)$ depends only on the number of balls other than ball $j$ in the boxes $x_j$ and $x'_j$. Thus, $f(x) - f(x^j) = f(x^i) - f(x^{ij})$. This proves that $G$ is an interaction graph for $f$.

Now, define $G'$ on $\mathcal{X}^{n+4}$ in exactly the same way as we defined $G$ on $\mathcal{X}^n$, that is, given $x \in \mathcal{X}^{n+4}$, $G'(x)$ puts an edge between $i$ and $j$ if and only if $x_i = x_j$. Again, it is trivial to check that $G'$ is symmetric and that $G'$ is an extension of $G$.

We now see that by the definition in Theorem 2.5, $\delta$ has the distribution of the number of balls in a typical box when we drop $n+4$ balls into $\alpha n$ boxes. Clearly, $\mathbb{E}(\delta^4) \le C\alpha^{-4}$ for some constant $C$ that does not depend on $n$. Finally, it is easy to check that $\sigma^2 \sim (\alpha e^{-1/\alpha} - (1+\alpha)e^{-2/\alpha})n$ as $n \to \infty$. The proof is now easy to complete using Theorem 2.5.  □

3.3. *Coverage processes.* Broadly speaking a stochastic coverage process is a random collection of (possibly overlapping) subsets of a metric space. The classic reference for the general theory of coverage processes is the book by Hall [28] (see also Chapter H in Aldous [1]).

We consider the following type of coverage process. Let $(\mathcal{X}, \rho)$ be a separable metric space endowed with a measure $\lambda$ (think of Euclidean space with Lebesgue measure) and suppose $X_1, \ldots, X_n$ are i.i.d. random points on $\mathcal{X}$ drawn according to some probability measure on $\mathcal{X}$. Fix some $\varepsilon > 0$ and let $\mathcal{R}$ be the random region covered by closed balls of radius $\varepsilon$ centered at $X_1, \ldots, X_n$ (our coverage process). Formally, if $\mathcal{B}(u, \varepsilon)$ denotes the closed ball of radius $\varepsilon$ centered at $u$, then

$$(3) \qquad \mathcal{R} = \bigcup_{i=1}^{n} \mathcal{B}(X_i, \varepsilon).$$

We will prove a general CLT for the area $\lambda(\mathcal{R})$. Of course, a large body of literature on this question already exists, but it is almost exclusively for processes on Euclidean spaces, where the analysis can be done by the big-blocks–small-blocks technique. The arguments are geometric in nature and do not extend to arbitrary metric spaces (e.g., manifolds). Moreover, the literature is silent on error bounds. For a discussion of the existing results and references, we refer to Section 3.4 of [28] (Theorem 3.5, in particular) and the notes at the end of Chapter 3 in the same book.

Here, we give a general normal approximation result with an error bound for the problem mentioned above. It comes as a very easy corollary of Theorem 2.5, possibly admitting extensions to more complex normal approximation problems in this area.

PROPOSITION 3.3. *Suppose we have $n$ i.i.d. points $X_1, \ldots, X_n$ on a separable metric space $(\mathcal{X}, \rho)$ endowed with a nonnegative Borel measure $\lambda$. Given $\varepsilon > 0$, define the set $\mathcal{R}$ as in (3). Put $M_\varepsilon = \sup_{u \in \mathcal{X}} \lambda(\mathcal{B}(u, \varepsilon))$ and $p_\varepsilon = \mathbb{P}(\rho(X_1, X_2) \le 2\varepsilon)$. Let $W = \lambda(\mathcal{R})$ and $\sigma_\varepsilon^2 = \mathrm{Var}(W)$. Then,*

$$\delta_W \le \frac{Cn^{1/2}M_\varepsilon^2(1 + np_\varepsilon)}{\sigma_\varepsilon^2} + \frac{nM_\varepsilon^3}{2\sigma_\varepsilon^3},$$

*where $\delta_W$ is the distance to normality defined in Definition 2.1 and $C$ is a universal constant.*

A bound like the above conveys no meaning unless applied to a concrete example. The simplest such example is the following. Let $\mathcal{X}$ be the unit square in $\mathbb{R}^2$ and $\varepsilon = n^{-1/2}$. Clearly, $M_\varepsilon \le C_1 n^{-1}$ and $p_\varepsilon \le C_2 n^{-1}$ for some constants $C_1$ and $C_2$ that do not depend on $n$. It can be shown (see [28], Theorem 3.4) that we also have $\sigma_\varepsilon^2 \ge C_3 n^{-1}$ for some positive constant $C_3$ free of $n$. Plugging these estimates into the above bound, we get $\delta_W \le Cn^{-1/2}$. Note that in this specific example, we may not get asymptotic normality if $\varepsilon$ decays faster than $n^{-1/2}$ as $n \to \infty$.

PROOF OF PROPOSITION 3.3. Given $x \in \mathcal{X}^n$, let $f(x) = \lambda(\mathcal{R}(x))$ and let $G(x)$ be the graph on $[n]$ that puts an edge between $i$ and $j$ if and only if $\rho(x_i, x_j) \le 2\varepsilon$. Let us verify that $G$ is an interaction rule for $f$.

Take any $x, x' \in \mathcal{X}^n$ and let $x^i$, $x^j$ and $x^{ij}$ be defined as in the beginning of Section 2.3. Let $N_j(x)$ be the set of neighbors of $x$ in the graph $G(x)$. Then, $f(x) - f(x^j) = \lambda(A) - \lambda(B)$, where

$$A = \mathcal{B}(x_j', \varepsilon) \Big\backslash \bigcup_{\ell \in N_j(x^j)} \mathcal{B}(x_\ell, \varepsilon) \quad \text{and} \quad B = \mathcal{B}(x_j, \varepsilon) \Big\backslash \bigcup_{\ell \in N_j(x)} \mathcal{B}(x_\ell, \varepsilon).$$

Now, if $\{i, j\}$ is not an edge in $G(x)$, $G(x^j)$, $G(x^i)$ and $G(x^{ij})$, then it is easy to see that $N_j(x) = N_j(x^i)$ and $N_j(x^j) = N_j(x^{ij})$. It follows that $f(x) = f(x^i)$ and $f(x^j) = f(x^{ij})$. Thus, $G$ is an interaction rule for $f$. The expression for $f(x) - f(x^j)$ also shows that $|f(x) - f(x^j)|$ is always bounded by the constant $M_\varepsilon$.

Next, given $x_1, \dots, x_{n+4} \in \mathcal{X}$, let $G'$ be defined in exactly the same way that $G$ was defined, that is, put include the edge $\{i, j\}$ if and only if $\rho(x_i, x_j) \le 2\varepsilon$. It is trivial to see that $G'$ is an extension of $G$ in the sense defined in Section 2.3. Thus, if $\delta$ is defined as in Theorem 2.5, then $\delta - 1 \sim \text{Binomial}(n+3, p_\varepsilon)$, where $p_\varepsilon = \mathbb{P}(\rho(X_1, X_2) \le 2\varepsilon)$. An application of Theorem 2.5 completes the proof. $\square$

3.4. *A CLT for nearest-neighbor statistics.* In a well-known 1983 paper, Bickel and Breiman [6] proved a central limit theorem for functionals of the form

$$(4) \qquad \frac{1}{\sqrt{n}} \sum_{\ell=1}^n (h(X_\ell, D_\ell) - \mathbb{E}h(X_\ell, D_\ell)),$$

where $X_1, \dots, X_n$ are i.i.d. random vectors following a probability density that is bounded and continuous on its support, $D_\ell := \min_{j \ne \ell} \|X_\ell - X_j\|$ is the distance between $X_\ell$ and its nearest neighbor and $h$ is a uniformly bounded and a.e. continuous function. Although the result looks very plausible, the proof is daunting. Indeed, as the authors put it, "*Our proof is long. We believe that this is due to the complexity of the problem.*" In short, their method can be described as a difficult multidimensional generalization of the familiar big-blocks–small-blocks method for analyzing $m$-dependent sequences.

Note that the existence of a density in the Bickel–Breiman theorem is a more restrictive assumption than it looks. For example, it precludes the possibility that the random variables are supported on some lower dimensional manifold, which may be quite important from a practical point of view.

In another widely cited work, Avram and Bertsimas [2] combined the Bickel–Breiman approach with the dependency graph technique of Baldi and Rinott [3] to yield CLTs for sums of edge lengths in various graphs arising from geometrical probability. A different method, originating from the work of Kesten and Lee [29], was used by Penrose and Yukich [32] to obtain a general CLT (with Kolmogorov distance error bound) for certain translation invariant functionals of uniformly distributed points and Poisson processes.

We have the following generalization of the Bickel–Breiman result, which, among other things, does away with the assumption that the $X_i$'s have a density with respect to Lebesgue measure. We also have an error bound, explicit up to a universal constant.

THEOREM 3.4. *Fix $n \geq 4$, $d \geq 1$, and $k \geq 1$. Suppose $X_1, \ldots, X_n$ are i.i.d. $\mathbb{R}^d$-valued random vectors with the property that $\|X_1 - X_2\|$ is a continuous random variable. Let $f : (\mathbb{R}^d)^n \to \mathbb{R}$ be a function of the form*

$$(5) \qquad f(x_1, \ldots, x_n) = \frac{1}{\sqrt{n}} \sum_{\ell=1}^n f_\ell(x_1, \ldots, x_n),$$

*where, for each $\ell$, $f_\ell(x_1, \ldots, x_n)$ is a function of only $x_\ell$ and its $k$ nearest neighbors. Suppose, for some $p \geq 8$, that $\gamma_p := \max_\ell \mathbb{E}|f_\ell(X_1, \ldots, X_n)|^p$ is finite. Let $W = f(X_1, \ldots, X_n)$ and $\sigma^2 = \mathrm{Var}(W)$. We then have the bound*

$$\delta_W \leq C \frac{\alpha(d)^3 k^4 \gamma_p^{2/p}}{\sigma^2 n^{(p-8)/2p}} + C \frac{\alpha(d)^3 k^3 \gamma_p^{3/p}}{\sigma^3 n^{(p-6)/2p}},$$

*where $\delta_W$ is the distance to normality defined in Definition 2.1, $\alpha(d)$ is the minimum number of $60°$ cones at the origin required to cover $\mathbb{R}^d$ and $C$ is a universal constant.*

REMARKS. (i) The assumption that the distribution of $\|X_1 - X_2\|$ does not have point masses is the bare minimal condition required to guarantee that the pairwise distances are all different (so that the nearest-neighbor orderings are uniquely defined). We believe that it is impossible to employ the big-blocks–small-blocks method under this minimal assumption, although it may be possible to formulate a version of the method that works when the $X_i$'s are supported on a sufficiently nice manifold.

(ii) The assumption concerning the $f_\ell$'s is also very weak. Unlike the Bickel–Breiman theorem, we do not require boundedness or continuity. Moreover, we do not even assume that the $f_\ell$'s are functions of only nearest-neighbor *distances*—they can be arbitrary functions of the nearest neighbors.

(iii) Like Theorem 2.2, the above result suffers from the deficiency that it does not say anything about $\sigma^2$. Again, as before, we think of that as a separate problem.

*Some applications.* (i) *Vertex degree in a geometric graph.* For a fixed $\varepsilon > 0$ and a given collection of points $x = (x_1, \ldots, x_n)$ in $\mathbb{R}^d$, the geometric graph $G(x, \varepsilon)$ is the graph on $x$ that puts edges between all pairs of vertices that are $\leq \varepsilon$ distance apart. Replacing $x$ by a collection $X = (X_1, \ldots, X_n)$ of i.i.d. random vectors, let $N_k$ be the number of points having vertex degree at least $k$ (where $k$ is fixed). This problem can be put in the context of Theorem 3.4 by defining $f_\ell(x) = 1$ if the distance between $x_\ell$ and its $k$th nearest neighbor is $\leq \varepsilon$, and $f_\ell(x) = 0$ otherwise. Then, $N_k = \sum_{\ell=1}^n f_\ell(x)$. Suppose all other terms are defined as in the statement of Theorem 3.4. Clearly, $\gamma_p \leq 1$ for all $p \geq 1$. Hence, we can take $p \to \infty$ and get

$$\delta_{N_k} \leq \frac{Ck^4}{\sigma^2 \sqrt{n}},$$

where $\sigma^2 = \mathrm{Var}(N_k)$, $\delta_{N_k}$ is the distance to normality defined in Definition 2.1 and $C$ is a constant depending on dimension $d$ and the distribution of the $X_\ell$'s. If $\varepsilon$ grows with $n$ at such a rate that $\sigma^2$ does not collapse to zero, then we get an $O(n^{-1/2})$ error bound for the Wasserstein distance. Incidentally, this example is quite well understood (see, e.g., Chapter 4 of [31]).

(ii) *Average nearest-neighbor distance.* Suppose $X_1, \ldots, X_n$ are i.i.d. random vectors in $\mathbb{R}^d$. Let $D_\ell$ be the distance of $X_\ell$ to its nearest-neighbor and $\bar{D} = \frac{1}{n} \sum_{\ell=1}^n D_\ell$ be the average nearest-neighbor distance. Assume that the support of the distribution of the $X_i$'s is $m$-dimensional, in the sense that the mass of $\varepsilon$-balls around any point is $\asymp \varepsilon^m$ as $\varepsilon \to 0$. Although a CLT for $\bar{D}$ could be proven using the Bickel–Breiman result if the $X_i$'s had a density with respect to Lebesgue measure, it does not work if we only assume that $\|X_1 - X_2\|$ has a continuous distribution.

Let $f_\ell = n^{1/m} D_\ell$ and $f = n^{-1/2} \sum_\ell f_\ell$. Then, for all $\varepsilon > 0$, we clearly have

$$\mathbb{P}(f_\ell(X) > \varepsilon) = \left(1 - \frac{C\varepsilon^m}{n}\right)^n \leq \exp(-C\varepsilon).$$

It follows that there is a constant $L \geq 1$ such that $\gamma_p^{1/p} \leq Lp$ for all $p \geq 1$. Along the same lines, it is not difficult to show that $\sigma^2 := \mathrm{Var}(f(X)) \asymp 1$ as

$n \to \infty$. Taking $p = \log n$, we get the bound

$$\delta_D \le \frac{C(\log n)^3}{\sqrt{n}},$$

where $\delta_D$ is the distance to normality defined in Definition 2.1 and $C$ is a constant depending on the dimension $d$ and the distribution of the $X_\ell$'s.

(iii) *The Levina–Bickel statistic.* In the preceding examples, we see that the error bound is effectively $O(n^{-1/2})$ when the summands have light tails. However, the $f_\ell$'s may be heavy-tailed in applications. A specific example of such a function is the recent "dimension estimator" of Levina and Bickel [30] which uses the distances to the first $k$ nearest neighbors to obtain an estimate of the so-called intrinsic dimension of a statistical data cloud. Explicitly, if $X_1, \ldots, X_n$ are i.i.d. random variables lying on a nice manifold of unknown dimension $m$ embedded in a higher-dimensional space $\mathbb{R}^d$, and $k$ is a positive integer $\ge 2$, then the Levina–Bickel estimate of $m$ with tuning parameter $k$ is given by the formula

$$(6) \qquad \hat{m}_k = \frac{1}{n} \sum_{\ell=1}^{n} \left( \frac{1}{k-1} \sum_{j=1}^{k-1} \log \frac{D_{\ell k}}{D_{\ell j}} \right)^{-1},$$

where $D_{\ell j}$ is the distance between $X_\ell$ and its $j$th nearest neighbor. In (6), we have $f_\ell(x) = (k-1)/g_\ell(x)$, where

$$g_\ell(x) = \sum_{j=1}^{k-1} \log \frac{D_{\ell k}(x)}{D_{\ell j}(x)}$$

and $D_{\ell j}(x)$ is the distance between $x_\ell$ and its $j$th nearest neighbor in the collection $x = (x_1, \ldots, x_n)$. It is argued in [30] that for large $n$, under appropriate assumptions, the distribution of $m \cdot g_\ell(X)$ can be approximated by the Gamma$(k, 1)$ distribution (recall that $m$ is the dimension of the manifold on which the data lie). It follows that

$$\mathbb{E}|f_\ell(X)|^{k-1} \le \frac{Cm^k(k-1)^k}{(k-1)!},$$

where $C$ is a constant that does not depend on $k$, $n$ and $m$. Putting $p = k-1$ in Theorem 3.4, we get

$$\delta_{\hat{m}_k} \le C \frac{\alpha(d)^3 k^3 m^2 (k\sigma + m)}{\sigma^3 n^{(k-9)/(2k-2)}},$$

where $\sigma^2 = \text{Var}(\sqrt{n}(\hat{m}_k - \mathbb{E}\hat{m}_k))$ and $\delta_{\hat{m}_k}$ is the distance to normality defined in Definition 2.1. Levina and Bickel ([30], Section 3) claim that for fixed $k$, they have a proof that $\sigma^2 \asymp 1$ as $n \to \infty$. This, combined with the above bound, implies a CLT for the Levina–Bickel statistic for $k > 9$.

PROOF OF THEOREM 3.4. For each $x = (x_1, \ldots, x_n) \in (\mathbb{R}^d)^n$, define a function $d_x$ on $[n] \times [n]$ as

$$(7) \qquad d_x(i, j) = \#\{\ell : \|x_i - x_\ell\| < \|x_i - x_j\|\}.$$

Our first task is to identify an interaction rule for functions of the form (5). Suppose $k$ is a fixed positive integer. Given any $x \in (\mathbb{R}^d)^n$, let $G(x)$ be the graph on $[n]$ that puts an edge between $i$ and $j$ if and only if there exists an $\ell$ such that $d_x(\ell, i) \leq k + 1$ and $d_x(\ell, j) \leq k + 1$. We claim that $G$ is a symmetric interaction rule for $f$.

To prove this claim, we begin with a simple observation: if $x, x' \in (\mathbb{R}^d)^n$ and $\ell, m \in [n]$ are such that $x_\ell = x'_\ell$ and $x_m = x'_m$, then

$$(8) \qquad |d_x(\ell, m) - d_{x'}(\ell, m)| \leq \#\{r : x_r \neq x'_r\}.$$

Now, fix some $x, x' \in (\mathbb{R}^d)^n$ and $i, j \in [n]$, where $i \neq j$. Define $x^i$, $x^j$ and $x^{ij}$ as in the definition of interaction between coordinates in Section 2.3. Suppose $\{i, j\}$ is not an edge in $G(x)$, $G(x^i)$, $G(x^j)$ and $G(x^{ij})$. We will show that for every $\ell$,

$$(9) \qquad f_\ell(x) - f_\ell(x^j) - f_\ell(x^i) + f_\ell(x^{ij}) = 0.$$

So, let us fix some $\ell \in [n]$. First, suppose that

$$(10) \qquad d_x(\ell, j) \leq k.$$

We claim that in this situation,

$$(11) \qquad f_\ell(x) = f_\ell(x^i) \quad \text{and} \quad f_\ell(x^j) = f_\ell(x^{ij}).$$

To show that, first note that since $\{i, j\} \notin G(x)$, we have

$$(12) \qquad d_x(\ell, i) > k + 1.$$

In particular, $i$ is different from $\ell$ and $j$. Thus, using (8) and (10), we see that $d_{x^i}(\ell, j) \leq k + 1$. Combining this with the hypothesis that $\{i, j\} \notin G(x^i)$, we get

$$(13) \qquad d_{x^i}(\ell, i) > k + 1.$$

From (12) and (13), it is easy to deduce that $x_\ell$ has the same set of $k$ nearest neighbors in both $x$ and $x^i$, hence that $f_\ell(x) = f_\ell(x^i)$.

Next, still assuming (10), suppose that $d_{x^j}(\ell, i) \leq k$. We show that this is impossible by considering two cases: (i) if $j = \ell$, this is clearly false because $\{i, j\} \notin G(x^j)$; (ii) if $j \neq \ell$, then by (8) and (12), we have $d_{x^j}(\ell, i) \geq k + 1$. Thus, under (10), we must have

$$(14) \qquad d_{x^j}(\ell, i) \geq k + 1.$$

Finally, still under (10), suppose we have $d_{x^{ij}}(\ell, i) \leq k$. Again, we show that this cannot be true under (10) by considering two cases: (i) if $\ell = j$, this cannot hold because $\{i, j\} \notin G(x^{ij})$; (ii) if $\ell \neq j$, then from (8) and (13), we get $d_{x^{ij}}(\ell, i) \geq k + 1$. Thus, under (10), we have

$$(15) \qquad\qquad d_{x^{ij}}(\ell, i) \geq k + 1.$$

From (14) and (15), it follows that $x_\ell$ has the same set of $k$ nearest neighbors in $x^j$ and $x^{ij}$. Therefore, $f_\ell(x^j) = f_\ell(x^{ij})$. This completes the proof of (11) under the hypothesis (10).

The symmetry in the problem now implies that (11) holds if $d_{x^j}(\ell, j) \leq k$ or $d_{x^i}(\ell, j) \leq k$, $d_{x^{ij}}(\ell, j) \leq k$. If none of these are true [i.e., $d_z(\ell, j) > k$ for $z = x, x^j, x^i, x^{ij}$], then we can directly deduce that the set of $k$ nearest neighbors of $x_\ell$ is the same in $x$ and $x^j$ and (separately) also in $x^i$ and $x^{ij}$, therefore

$$f_\ell(x) = f_\ell(x^j) \quad \text{and} \quad f_\ell(x^i) = f_\ell(x^{ij}).$$

Combining the cases, the proof of (9) is now complete.

Thus, we have proven the claim that $G$ is an interaction rule for $f$. Clearly, $G$ is symmetric. A symmetric extension of $G$ to $(\mathbb{R}^d)^{n+4}$ is easily constructed as follows. Given any vector $x \in (\mathbb{R}^d)^{n+4}$, let $G'(x)$ be the graph on $[n + 4]$ that puts an edge between $i$ and $j$ if and only if there exists an $\ell \in [n + 4]$ such that $d_x(\ell, i) \leq k + 5$ and $d_x(\ell, j) \leq k + 5$. To see this, note that if $\{i, j\} \in G(x_1, \ldots, x_n)$, then there exists some $\ell$ such that $x_i$ and $x_j$ are both among the $k + 1$ nearest neighbors of $x_\ell$ in the set $\{x_1, \ldots, x_n\}$. After the addition of four more points to this set, $x_i$ and $x_j$ will still be members of the set of $k + 5$ nearest neighbors of $x_\ell$. This proves that $G'$ is an extension of $G$, and it is obviously symmetric.

Now, for every $x \in \mathbb{R}^d$ and $1 \leq j \leq n$, let

$$N_j(x) := \{\ell : d_x(\ell, j) \leq k\}.$$

As we have noted before, if $\ell \notin N_j(x) \cup N_j(x^j)$, then $x_\ell$ has the same set of $k$ nearest neighbors in both $x$ and $x^j$, therefore $f_\ell(x) = f_\ell(x^j)$. Thus,

$$f(x) - f(x^j) = \sum_{\ell \in N_j(x) \cup N_j(x^j)} n^{-1/2}(f_\ell(x) - f_\ell(x^j)).$$

It follows from standard geometrical arguments (see, e.g., [42], page 102) and the assumption that $\|X_1 - X_2\|$ is a continuous r.v. that $|N_j(x) \cup N_j(x^j)| \leq 2\alpha(d)k$, irrespective of $n$ and $x$, where $\alpha(d)$ is the minimum number of $60°$ cones at the origin required to cover $\mathbb{R}^d$. Thus, if we let

$$M_f := \max_\ell |f_\ell(X)| \vee \max_{j, \ell} |f_\ell(X^j)|,$$

then the random variable $M$ in the statement of Theorem 2.5 can be bounded by $4n^{-1/2}\alpha(d)kM_f$ in this problem. Next, note that for any $p \geq 8$,

$$\mathbb{E}(M_f^8) \leq [\mathbb{E}(M_f^p)]^{8/p}$$

$$\leq \left[ \sum_{\ell} \mathbb{E}|f_{\ell}(X)|^p + \sum_{j,\ell} \mathbb{E}|f_{\ell}(X^j)|^p \right]^{8/p} \leq (n^2 + n)^{8/p}\gamma_p^{8/p}.$$

Similarly, one can show that $\mathbb{E}|\Delta_j f(X)|^3 \leq C\alpha(d)^3 k^3 n^{-3/2}(n\gamma_p)^{3/p}$. Finally, note that by the same geometrical observation as mentioned before, the maximum degree of $G'(X)$ is bounded by $\alpha(d)(k+1)(k+5)$. The proof is now completed by combining the bounds for all the terms and using Theorem 2.5. $\square$

## 4. Proofs of the main results.

4.1. *Proof of Theorem 2.2.* Let us begin with the observation that, without loss of generality, we can replace $f$ by $\sigma^{-1}f$ and then assume that $\sigma^2 = 1$. Henceforth, we will work under that assumption. The argument is divided into a sequence of lemmas. Lemmas 2.3 and 2.4 (already stated in Section 2) and Lemma 4.1 are original contributions of this paper, while Lemma 4.2 goes back to Stein [41].

PROOF OF LEMMA 2.3. Consider the sum

$$\sum_{A \subsetneq [n]} \frac{1}{\binom{n}{|A|}(n - |A|)} \sum_{j \notin A} \Delta_j f(X^A)$$

$$= \sum_{A \subsetneq [n]} \frac{1}{\binom{n}{|A|}(n - |A|)} \sum_{j \notin A} (f(X^A) - f(X^{A \cup j})).$$

Clearly, this is a linear combination of $\{f(X^A), A \subseteq [n]\}$. It is a matter of simple verification that the positive and negative coefficients of $f(X^A)$ in this linear combination cancel out except when $A = [n]$ or $A = \varnothing$. In fact, the above expression is identically equal to $f(X) - f(X')$.

Now, fix $A$ and $j \notin A$, and let $U = g(X)\Delta_j f(X^A)$. $U$ is then a function of the random vectors $X$ and $X'$. The joint distribution of $(X, X')$ remains unchanged if we interchange $X_j$ and $X_j'$. Under this operation, $U$ changes to $U' := -g(X^j)\Delta_j f(X^A)$. Thus,

$$\mathbb{E}(U) = \mathbb{E}(U') = \tfrac{1}{2}\mathbb{E}(U + U') = \tfrac{1}{2}\mathbb{E}[\Delta_j g(X)\Delta_j f(X^A)].$$

Combining these observations, we get

$$\mathrm{Cov}(g(X), f(X)) = \mathbb{E}[g(X)(f(X) - f(X'))]$$

$$= \sum_{A \subsetneq [n]} \frac{1}{\binom{n}{|A|}(n-|A|)} \sum_{j \notin A} \mathbb{E}[g(X)\Delta_j f(X^A)]$$

$$= \frac{1}{2} \sum_{A \subsetneq [n]} \frac{1}{\binom{n}{|A|}(n-|A|)} \sum_{j \notin A} \mathbb{E}[\Delta_j g(X)\Delta_j f(X^A)].$$

This completes the proof of the lemma. □

PROOF OF LEMMA 2.4. For each $A \subseteq [n]$ and $j \notin A$, let

$$R_{A,j} := \Delta_j(\varphi \circ f)(X)\Delta_j f(X^A)$$

and

$$\tilde{R}_{A,j} := \varphi'(f(X))\Delta_j f(X)\Delta_j f(X^A).$$

By Lemma 2.3 with $g = \varphi \circ f$, we have

$$(16) \qquad \mathbb{E}(\varphi(W)W) = \frac{1}{2} \sum_{A \subsetneq [n]} \frac{1}{\binom{n}{|A|}(n-|A|)} \sum_{j \notin A} \mathbb{E}(R_{A,j}).$$

By the mean value theorem, we have

$$\mathbb{E}|R_{A,j} - \tilde{R}_{A,j}| \le \frac{\|\varphi''\|_\infty}{2} \mathbb{E}|(\Delta_j f(X))^2 \Delta_j f(X^A)|$$

$$(17)$$

$$\le \frac{\|\varphi''\|_\infty}{2} \mathbb{E}|\Delta_j f(X)|^3 \qquad \text{(by Hölder's inequality)}.$$

Now, from the definition of $T$, we have

$$(18) \qquad \varphi'(W)T = \frac{1}{2} \sum_{A \subsetneq [n]} \frac{1}{\binom{n}{|A|}(n-|A|)} \sum_{j \notin A} \tilde{R}_{A,j}.$$

Combining (16), (17) and (18), we get

$$|\mathbb{E}(\varphi(W)W) - \mathbb{E}(\varphi'(W)T)|$$

$$= \left| \frac{1}{2} \sum_{A \subsetneq [n]} \frac{1}{\binom{n}{|A|}(n-|A|)} \sum_{j \notin A} \mathbb{E}(R_{A,j} - \tilde{R}_{A,j}) \right|$$

$$\le \frac{\|\varphi''\|_\infty}{4} \sum_{A \subsetneq [n]} \frac{1}{\binom{n}{|A|}(n-|A|)} \sum_{j \notin A} \mathbb{E}|\Delta_j f(X)|^3$$

$$= \frac{\|\varphi''\|_\infty}{4} \sum_{j=1}^{n} \mathbb{E}|\Delta_j f(X)|^3.$$

This completes the proof of the lemma. □

LEMMA 4.1.   *Let $W$ be as above. For any $\varphi \in C^2(\mathbb{R})$ such that $\|\varphi'\|_\infty \leq 1$ and $\|\varphi''\|_\infty \leq 2$, we have*

$$|\mathbb{E}(\varphi(W)W) - \mathbb{E}(\varphi'(W))| \leq [\mathrm{Var}(\mathbb{E}(T|W))]^{1/2} + \frac{1}{2}\sum_{j=1}^n \mathbb{E}|\Delta_j f(X)|^3.$$

PROOF.   Note that by putting $g = f$ in Lemma 2.3, we get $\mathbb{E}(T) = \mathbb{E}(W^2) = 1$. Since $\|\varphi'\|_\infty \leq 1$, this gives

$$|\mathbb{E}(\varphi'(W)T) - \mathbb{E}(\varphi'(W))| \leq \mathbb{E}|\mathbb{E}(T|W) - 1| \leq [\mathrm{Var}(\mathbb{E}(T|W))]^{1/2}.$$

The proof is completed by applying Lemma 2.4.   $\square$

LEMMA 4.2.   *Suppose $h\colon \mathbb{R} \to \mathbb{R}$ is an absolutely continuous function with bounded derivative. Let $Z \sim N(0,1)$. There then exists a solution to the differential equation*

$$\varphi'(x) - x\varphi(x) = h(x) - \mathbb{E}h(Z)$$

*that satisfies $\|\varphi'\|_\infty \leq \sqrt{\frac{2}{\pi}}\|h'\|_\infty$ and $\|\varphi''\|_\infty \leq 2\|h'\|_\infty$.*

REMARK.   It is not difficult to show that both constants are sharp. For a different proof of the bound on $\|\varphi'\|_\infty$, see Lemma 1 in [33]. The bound on $\|\varphi''\|_\infty$ is due to Stein ([41], page 27). Easier proofs with suboptimal constants can be found in Chen and Shao ([14], Chapter 1, Lemma 2.3).

PROOF OF LEMMA 4.2.   It can be verified that the function

$$\varphi(x) = e^{x^2/2}\int_{-\infty}^x e^{-t^2/2}(h(t) - \mathbb{E}h(Z))\,dt$$

$$= -e^{x^2/2}\int_x^\infty e^{-t^2/2}(h(t) - \mathbb{E}h(Z))\,dt$$

is a solution. Stein ([41], page 25, Lemma 3) proves that $\|\varphi''\|_\infty \leq 2\|h'\|_\infty$. The inequality $\|\varphi'\|_\infty \leq \|h'\|_\infty$ can also be derived using Stein's proof of the other inequality. We carry out the steps below. First, it is easy to verify that

$$h(x) - \mathbb{E}h(Z) = \int_{-\infty}^x h'(z)\Phi(z)\,dz - \int_x^\infty h'(z)(1 - \Phi(z))\,dz,$$

where $\Phi$ is the standard Gaussian c.d.f. Again, as proven in Stein ([41], page 27),

$$\varphi(x) = -\sqrt{2\pi}e^{x^2/2}(1 - \Phi(x))\int_{-\infty}^x h'(z)\Phi(z)\,dz$$

$$- \sqrt{2\pi}e^{x^2/2}\Phi(x)\int_x^\infty h'(z)(1 - \Phi(z))\,dz.$$

Combining, we see that

$$\varphi'(x) = x\varphi(x) + h(x) - \mathbb{E}h(Z)$$
$$= (1 - \sqrt{2\pi}xe^{x^2/2}(1 - \Phi(x)))\int_{-\infty}^{x} h'(z)\Phi(z)\,dz$$
$$- (1 + \sqrt{2\pi}xe^{x^2/2}\Phi(x))\int_{x}^{\infty} h'(z)(1 - \Phi(z))\,dz.$$

It follows that

$$\|\varphi'\|_\infty \leq \|h'\|_\infty \sup_{x \in \mathbb{R}} \bigg( |1 - \sqrt{2\pi}xe^{x^2/2}(1 - \Phi(x))| \int_{-\infty}^{x} \Phi(z)\,dz$$
$$+ |1 + \sqrt{2\pi}xe^{x^2/2}\Phi(x)| \int_{x}^{\infty} (1 - \Phi(z))\,dz \bigg).$$

Using integration by parts, we get

$$\int_{-\infty}^{x} \Phi(z)\,dz = x\Phi(x) + \frac{e^{-x^2/2}}{\sqrt{2\pi}}$$

and

$$\int_{x}^{\infty} (1 - \Phi(z))\,dz = -x(1 - \Phi(x)) + \frac{e^{-x^2/2}}{\sqrt{2\pi}}.$$

Thus, we have

$$\|\varphi'\|_\infty \leq \|h'\|_\infty \sup_{x \in \mathbb{R}} \bigg( |1 - \sqrt{2\pi}xe^{x^2/2}(1 - \Phi(x))| \Big( x\Phi(x) + \frac{e^{-x^2/2}}{\sqrt{2\pi}} \Big)$$
$$+ |1 + \sqrt{2\pi}xe^{x^2/2}\Phi(x)| \Big( -x(1 - \Phi(x)) + \frac{e^{-x^2/2}}{\sqrt{2\pi}} \Big) \bigg).$$

It is a calculus exercise to verify that the term inside the brackets attains its maximum at $x = 0$, where its value is $\sqrt{2/\pi}$.   □

PROOF OF THEOREM 2.2.   Take any $h$ with $\|h'\|_\infty \leq 1$. Let $\varphi$ be a solution to $\varphi'(x) - x\varphi(x) = h(x) - \mathbb{E}h(Z)$. Then,

$$\mathbb{E}h(W) - \mathbb{E}h(Z) = \mathbb{E}(\varphi'(W)) - \mathbb{E}(W\varphi(W)).$$

By Lemma 4.2, $\|\varphi'\|_\infty \leq \sqrt{\frac{2}{\pi}} \leq 1$ and $\|\varphi''\|_\infty \leq 2$. The proof is completed by applying Lemma 4.1.   □

4.2. *Proof of Theorem 2.5.* By Theorem 2.2, our task reduces to obtaining a bound on $\operatorname{Var}(\mathbb{E}(T|X))$, where $T$ is defined in (1). However, the situation in Theorem 2.5 is too complex to admit a direct computation of the variance. To circumvent this problem, we will use the following well-known martingale bound for the variance of an arbitrary function of independent random variables. This is known as the *Efron–Stein inequality* in the statistics literature.

LEMMA 4.3 ([17, 39]). *Let $Z = g(Y_1, \ldots, Y_m)$ be a function of independent random objects $Y_1, \ldots, Y_m$. Let $Y_i'$ be an independent copy of $Y_i$, $i = 1, \ldots, m$. Then,*

$$\operatorname{Var}(Z) \leq \frac{1}{2} \sum_{i=1}^{m} \mathbb{E}[(g(Y_1, \ldots, Y_{i-1}, Y_i', Y_{i+1}, \ldots, Y_m) - g(Y_1, \ldots, Y_m))^2].$$

We will combine this inequality with another simple inequality that we were unable to locate in the literature.

LEMMA 4.4. *If $X$ and $X'$ are independent random objects, then for any square integrable function $U = g(X, X')$, we have the inequality*

$$\operatorname{Var}(\mathbb{E}(U|X)) \leq \mathbb{E}(\operatorname{Var}(U|X')).$$

PROOF. The proof is based on a simple application of Jensen's inequality. We just note that by the independence of $X$ and $X'$, we have $\mathbb{E}(\mathbb{E}(U|X')|X) = \mathbb{E}(U)$ and therefore

$$\begin{aligned}
\operatorname{Var}(\mathbb{E}(U|X)) &= \mathbb{E}(\mathbb{E}(U|X) - \mathbb{E}(U))^2 \\
&= \mathbb{E}(\mathbb{E}(U - \mathbb{E}(U|X')|X))^2 \\
&\leq \mathbb{E}(U - \mathbb{E}(U|X'))^2 = \mathbb{E}(\operatorname{Var}(U|X')).
\end{aligned}$$

This completes the proof of the lemma. □

Now, recall the definitions of $\Delta_j$, $T_A$ and $T$ from Section 2, and the normal approximation bound in terms of $\operatorname{Var}(\mathbb{E}(T|X))$ in Theorem 2.2. We will prove the following upper bound on $\operatorname{Var}(\mathbb{E}(T_A|X))$.

LEMMA 4.5. *With everything defined as before, we have*

$$\operatorname{Var}(\mathbb{E}(T_A|X)) \leq C\mathbb{E}(M^8)^{1/2}\mathbb{E}(\delta^4)^{1/2}\sqrt{n(n-|A|)},$$

*where $C$ is a universal constant.*

This is a good place to declare the convention that throughout the remainder of this section, $C$ will denote numerical constants that do not depend on anything else and the value of $C$ may change from line to line.

Before proving Lemma 4.5, we need to finish an important task.

PROOF OF THEOREM 2.2.   Lemma 4.5, combined with Theorem 2.2, completes the proof of Theorem 2.5 as follows. First, note that by the definition (1) of $T$ and Minkowski's inequality, we have

$$[\mathrm{Var}(\mathbb{E}(T|X))]^{1/2} \le \frac{1}{2} \sum_{A \subsetneq [n]} \frac{[\mathrm{Var}(\mathbb{E}(T_A|X))]^{1/2}}{\binom{n}{|A|}(n-|A|)}.$$

Substituting the bound from Lemma 4.5 into the above expression, we get

$$[\mathrm{Var}(\mathbb{E}(T|X))]^{1/2} \le C\mathbb{E}(M^8)^{1/2}\mathbb{E}(\delta^4)^{1/2} \sum_{A \subsetneq [n]} \frac{n^{1/4}(n-|A|)^{1/4}}{\binom{n}{|A|}(n-|A|)}$$

$$= C\mathbb{E}(M^8)^{1/2}\mathbb{E}(\delta^4)^{1/2} \sum_{k=1}^{n} n^{1/4} k^{-3/4}$$

$$\le C\mathbb{E}(M^8)^{1/2}\mathbb{E}(\delta^4)^{1/2} n^{1/2}.$$

This completes the proof of Theorem 2.5.   □

Our main job now is to prove Lemma 4.5. Let us begin with a simple lemma about symmetric graphical rules.

LEMMA 4.6.   *Suppose $G$ is a symmetric graphical rule on $\mathcal{X}^n$ and $X = (X_1, X_2, \ldots, X_n)$ is a vector of independent and identically distributed $\mathcal{X}$-valued random variables. Let $d_1$ be the degree of the vertex 1 in $G(X)$. Take any $k \le n-1$ and let $i, i_1, i_2, \ldots, i_k$ be any collection of $k+1$ distinct elements of $[n]$. Then,*

(19)          $$\mathbb{P}(\{i, i_\ell\} \in G(X) \text{ for each } 1 \le \ell \le k) = \frac{\mathbb{E}((d_1)_k)}{(n-1)_k},$$

*where $(r)_k$ stands for the product $r(r-1)\cdots(r-k+1)$.*

PROOF.   Since $G$ is a symmetric rule and the $X_i$'s are i.i.d., the quantity

$$\mathbb{P}(\{i, i_\ell\} \in G(X) \text{ for all } 1 \le \ell \le k)$$

does not depend on the specific choice of $i, i_1, \ldots, i_k$. Hence,

$$\mathbb{P}(\{i, i_\ell\} \in G(X) \text{ for each } 1 \le \ell \le k)$$

$$= \frac{1}{(n-1)_k} \sum \mathbb{P}(\{i, j_\ell\} \in G(X) \text{ for each } 1 \le \ell \le k),$$

where the sum is taken over all choices of distinct $j_1, \ldots, j_k$ in $[n] \setminus \{i\}$. Finally, note that

$$\sum \mathbb{I}\{\{i, j_\ell\} \in G(X) \text{ for each } 1 \leq \ell \leq k\} = (d_i)_k,$$

where $d_i$ is the degree of the vertex $i$. Again, by symmetry, $d_i$ and $d_1$ have the same distribution. This completes the argument. $\square$

Proof of Lemma 4.5. Fix a set $A \subsetneq [n]$. For each $j \notin A$, let

$$R_j = \Delta_j f(X) \Delta_j f(X^A)$$
$$= (f(X) - f(X^j))(f(X^A) - f(X^{A \cup j})).$$

Now, let $Y = (Y_1, \ldots, Y_n)$ be another copy of $X$, independent of both $X$ and $X'$. Fix $1 \leq i \leq n$. Let

$$\tilde{X} = (X_1, \ldots, X_{i-1}, Y_i, X_{i+1}, \ldots, X_n).$$

Similarly, for each $B \subseteq [n]$, define $\tilde{X}^B$ by replacing $X_i$ with $Y_i$ in $X^B$. Explicitly, if $i \notin B$, then

$$\tilde{X}^B = (X_1^B, \ldots, X_{i-1}^B, Y_i, X_{i+1}^B, \ldots, X_n^B),$$

whereas if $i \in B$, then $\tilde{X}^B = X^B$. With this notation, let

$$R_{ji} = (f(\tilde{X}) - f(\tilde{X}^j))(f(\tilde{X}^A) - f(\tilde{X}^{A \cup j})),$$

and put

$$h_i := \mathbb{E}\left(\sum_{j \notin A} (R_j - R_{ji})\right)^2.$$

It follows from a combination of Lemmas 4.3 and 4.4 that

$$(20) \qquad \text{Var}(\mathbb{E}(T_A|X)) \leq \mathbb{E}(\text{Var}(T_A|X')) \leq \frac{1}{2} \sum_{i=1}^{n} h_i.$$

Let us now proceed to bound $h_i$. First, take some $j \notin A \cup i$ and let

$$d_{ji}^1 = \mathbb{I}\{\{i, j\} \in G(X)\},$$
$$d_{ji}^2 = \mathbb{I}\{\{i, j\} \in G(X^j)\},$$
$$d_{ji}^3 = \mathbb{I}\{\{i, j\} \in G(\tilde{X})\}$$

and

$$d_{ji}^4 = \mathbb{I}\{\{i, j\} \in G(\tilde{X}^j)\}.$$

Now, suppose that for a particular realization, we have $d_{ji}^1 = d_{ji}^2 = d_{ji}^3 = d_{ji}^4 = 0$. Since $G$ is an interaction rule for $f$, this event implies that

$$f(X) - f(X^j) = f(\tilde{X}) - f(\tilde{X}^j).$$

If we now take $X^A$ instead of $X$ and $\tilde{X}^A$ instead of $\tilde{X}$, and define $e_{ji}^1$, $e_{ji}^2$, $e_{ji}^3$ and $e_{ji}^4$ analogously, then the event $e_{ji}^1 = e_{ji}^2 = e_{ji}^3 = e_{ji}^4 = 0$ implies that

$$f(X^A) - f(X^{A \cup j}) = f(\tilde{X}^A) - f(\tilde{X}^{A \cup j}),$$

irrespective of whether or not $i \in A$. Now, let

$$L_i := \max_{j \notin A} |\Delta_j f(X) \Delta_j f(X^A) - \Delta_j f(\tilde{X}) \Delta_j f(\tilde{X}^A)|.$$

From the preceding observations, we see that for $j \notin A \cup i$,

$$|R_j - R_{ji}| \le L_i \sum_{k=1}^{4} (d_{ji}^k + e_{ji}^k).$$

When $i \notin A$ and $j = i$, we simply have $|R_j - R_{ji}| \le L_i$. Applying the Cauchy–Schwarz inequality, we now get

$$(21) \qquad h_i \le \left[ \mathbb{E}(L_i^4) \mathbb{E}\left( \mathbb{I}\{i \notin A\} + \sum_{j \notin A \cup i} \sum_{k=1}^{4} (d_{ji}^k + e_{ji}^k) \right)^4 \right]^{1/2}.$$

Now, by the inequality $(\sum_{i=1}^{r} a_i)^4 \le r^3 \sum_{i=1}^{r} a_i^4$, we have

$$\mathbb{E}\left( \mathbb{I}\{i \notin A\} + \sum_{j \notin A \cup i} \sum_{k=1}^{4} (d_{ji}^k + e_{ji}^k) \right)^4$$

$$\le 9^3 \mathbb{I}\{i \in A\} + 9^3 \sum_{k=1}^{4} \mathbb{E}\left( \sum_{j \notin A \cup i} d_{ji}^k \right)^4 + 9^3 \sum_{k=1}^{4} \mathbb{E}\left( \sum_{j \notin A \cup i} e_{ji}^k \right)^4.$$

To get a bound for the above terms, first consider the $d^1$ term. It follows directly from Lemma 4.6 that for any $j, k, l$ and $m$,

$$\mathbb{E}(d_{ji}^1 d_{ki}^1 d_{li}^1 d_{mi}^1) \le C \frac{\mathbb{E}(\delta_1^r)}{n^r},$$

where $r =$ the number of distinct indices among $j, k, l, m$ and $\delta_1$ is the degree of the vertex 1 in $G(X)$. Recall the definition of $\delta$ from the statement of the theorem and observe that $\delta \ge \delta_1 + 1$. It is now easy to deduce that

$$\mathbb{E}\left( \sum_{j \notin A \cup i} d_{ji}^1 \right)^4 \le C \mathbb{E}(\delta^4) \frac{n - |A|}{n}.$$

Now, consider the problem of bounding $\mathbb{E}(d_{ji}^2 d_{ki}^2 d_{li}^2 d_{mi}^2)$. First, suppose $j$, $k$, $l$ and $m$ are distinct. Let $\tilde{X}$ be the random vector on $\mathcal{X}^{n+4}$ defined as

$$\tilde{X} := (X_1, \ldots, X_n, X_j', X_k', X_l', X_m').$$

Note that if $d_{ji}^2 = d_{ki}^2 = d_{li}^2 = d_{mi}^2 = 1$, then $\{i, n+1\}$, $\{i, n+2\}$, $\{i, n+3\}$ and $\{i, n+4\}$ are all edges in the extended graph $G'(\tilde{X})$. Since $G'$ is a symmetric rule and the components of $\tilde{X}$ are i.i.d., it again follows from Lemma 4.6 that

$$\mathbb{E}(d_{ji}^2 d_{ki}^2 d_{li}^2 d_{mi}^2) \le C \frac{\mathbb{E}(\delta^4)}{n^4}.$$

Now, suppose $j, k, l$ are distinct, but $m = l$. Let $s$ be an element of $[n]$ different from $j$, $k$ and $l$. Define

$$\tilde{X} := (X_1, \ldots, X_n, X_j', X_k', X_l', X_s')$$

and proceed as before to conclude that, in this case,

$$\mathbb{E}(d_{ji}^2 d_{ki}^2 d_{li}^2 d_{mi}^2) = \mathbb{E}(d_{ji}^2 d_{ki}^2 d_{li}^2) \le C \frac{\mathbb{E}(\delta^3)}{n^3}.$$

In general, if $r$ is the number of distinct elements among $j, k, l, m$, then

$$\mathbb{E}(d_{ji}^2 d_{ki}^2 d_{li}^2 d_{mi}^2) \le C \frac{\mathbb{E}(\delta^r)}{n^r}.$$

From this, we get

$$\mathbb{E}\left(\sum_{j \notin A \cup i} d_{ji}^2\right)^4 \le C \mathbb{E}(\delta^4) \frac{n - |A|}{n}.$$

The $d^3, e^1$ and $e^3$ terms can be given the same bound as the $d^1$ term, while the $d^4$, $e^2$ and $e^4$ terms are similar to the $d^2$ term. Combining, we get

$$\mathbb{E}\left(\mathbb{I}\{i \notin A\} + \sum_{j \notin A} \sum_{k=1}^4 (d_j^k + e_j^k)\right)^4 \le C \mathbb{E}(\delta^4)\left(\mathbb{I}\{i \notin A\} + \frac{n - |A|}{n}\right).$$

It is easy to show, using the Cauchy–Schwarz inequality, that $\mathbb{E}(L_i^4) \le C \mathbb{E}(M^8)$, where $M = \max_j |\Delta_j f(X)|$. Using these bounds in (21) and the inequality $\sqrt{x + y} \le \sqrt{x} + \sqrt{y}$, we get

$$h_i \le C \mathbb{E}(M^8)^{1/2} \mathbb{E}(\delta^4)^{1/2}\left(\mathbb{I}\{i \notin A\} + \sqrt{\frac{n - |A|}{n}}\right).$$

Substituting this bound in (20), we get

$$\mathrm{Var}(\mathbb{E}(T_A|X)) \le C \mathbb{E}(M^8)^{1/2} \mathbb{E}(\delta^4)^{1/2}(n - |A| + \sqrt{n(n - |A|)})$$

$$\le C \mathbb{E}(M^8)^{1/2} \mathbb{E}(\delta^4)^{1/2} \sqrt{n(n - |A|)}.$$

This completes the proof.  $\square$

**Acknowledgments.** The author wishes to thank Persi Diaconis, Peter Bickel, Yuval Peres and the anonymous referee for many useful comments and suggestions.

Department of Statistics
University of California at Berkeley
367 Evans Hall #3860
Berkeley, California 94720-3860
USA
E-mail: sourav@stat.berkeley.edu
URL: http://www.stat.berkeley.edu/˜sourav/